\documentclass[12pt]{amsproc}
\usepackage[utf8]{inputenc}
\usepackage[total={6.5in,9in},top=1in, left=1in, right=1in, bottom=1in]{geometry}
\usepackage{fancyhdr}
\usepackage{libertine}
\usepackage[libertine]{newtxmath}
\usepackage{color,amsmath,mathtools,graphicx}
\usepackage{bbold}
\usepackage{epsfig,fancyhdr}
\usepackage{dsfont} 
\usepackage{graphicx,amsmath,mathtools,float}
\usepackage{epsfig,color,fancyhdr,setspace}
\usepackage{dsfont}
\usepackage{epstopdf}
\usepackage{import}
\usepackage{lineno}
\usepackage{stackrel,dsfont}
\usepackage[allcolors = blue,colorlinks]{hyperref}
\usepackage{hyperref}
\usepackage[abs]{overpic}
\usepackage[center]{caption}
\usepackage{subcaption}
\usepackage{tikz-cd}
\usepackage{enumitem}
\usepackage[utf8]{inputenc}
\usepackage[T1]{fontenc}

\newtheorem{theorem}{Theorem}[section]

\newtheorem{definition}[theorem]{Definition}

\newtheorem{remark}[theorem]{Remark}

\newtheorem{conjecture}[theorem]{Conjecture}

\title{A Counterexample to the 
Generalisation of Witten's Conjecture}

 \dedicatory{Dedicated to the memory of Zbigniew Oziewicz}

\author{Rhea Palak Bakshi}
\address{Institute for Theoretical Studies, ETH Zürich, Zürich, Switzerland}
\email{rheapalak.bakshi@eth-its.ethz.ch}
\thanks{The author would like to thank Hanna Makaruk and Robert Owczarek for organising several AMS Special Sessions on Inverse Problems over the years and for their endless support. It was at a session they had organised in Portland in 2018 that the author first had the opportunity to meet Zbigniew Oziewicz.}

\keywords{Knot, 3-manifold, link invariant, $3$-manifold invariant, Kauffman bracket skein module.}
\subjclass[2020]{Primary: 57K31. Secondary: 57K10}

\begin{document}

\begin{abstract}

This note provides a counterexample to a conjecture by Marché about the structure of the Kauffman bracket skein module for closed compact oriented $3$-manifolds over the ring of Laurent polynomials. 

\end{abstract}

\maketitle

\section{Introduction}\label{intromarche}

Classical knot and link theory is the study of embeddings of circles in the $3$-sphere up to ambient isotopy. Much of this theory concerns itself with determining which knots and links are equivalent and which are not. Over the years several link invariants have been discovered, which have been reasonably successful in telling knots and links apart from one another. Link polynomials, such as the Alexander-Conway polynomial \cite{Con}, Jones polynomial \cite{jonesannals}, HOMFLYPT polynomial \cite{homfly,homflypt}, Kauffman $2$-variable polynomial \cite{2variablekauffman}, Kauffman bracket polynomial \cite{sm&j}, and Dubrovnik polynomial \cite{2variablekauffman}, form an important class of these invariants. What these polynomials have in common is that their definitions are built upon skein relations. A skein relation is a linear relation between two or more link diagrams that are exactly the same everywhere except in the neighbourhood of a crossing where they differ. Skein modules arise naturally when one tries to find analogues of all these polynomial link invariants for links in any arbitrary $3$-manifold. \\

Given an oriented\footnote{In order to replace $S^3$ with an arbitrary $3$-manifold $M$ we consider links which are identical outside a ball $B^3$ in $M$ and which appear inside $B^3$ as various prescribed tangles. To be able to distinguish one link from another locally inside $B^3$ we usually assume that $M$ is oriented.} $3$-manifold $M$ and a commutative ring with unity, there are several skein modules that can be associated to $M$, each of which tries to capture some information about the knot theory that $M$ admits. One reason for this plethora of definitions of skein modules is that we can consider links in $M$ to be either oriented or unoriented, framed or unframed, or up to isotopy, homotopy, or homology, for example. What all these definitions have in common, however, is that they generalize the skein theory of all the various polynomial link invariants for links in $S^3$ to arbitrary $3$-manifolds. Thus, the theory of skein modules, as introduced by Przytycki in \cite{smof3}, can be seen as the natural extension of the Alexander-Conway, HOMFLYPT, Kauffman bracket, and Jones polynomial link invariants in $S^3$ to links in arbitrary $3$-manifolds. In fact, we allow skein relations involving arbitrary linear combinations of tangles in an oriented $3$-manifold. For a brief survey of skein modules see \cite{ch63ency}. \\ 

The skein module based on the Kauffman bracket polynomial is by far the most extensively studied skein module of all. It has strong ties with algebraic and hyperbolic geometry via $SL(2, \mathbb{C})$ character varieties \cite{sl2cbullock}, quantum Teichmüller spaces \cite{bw1,fokchekhov,lequantum,kashaev}, and quantum cluster algebras \cite{mullerskein}. It is also instrumental in the study of the $AJ$ conjecture, which relates the coloured Jones polynomial and the $A$-polynomial of a knot \cite{FGL,garou,2bk}. Moreover, it is used in the combinatorial construction \cite{Li4,BHMV1} of the Witten-Reshetikhin-Turaev Topological Quantum Field Theories. Due to its fundamental role in this construction and its intricate ties to knot theory, geometric topology, and hyperbolic geometry, the Kauffman bracket skein module has become pivotal in the study of quantum topology.\\

\begin{definition}

Let $M$ be an oriented $3$-manifold, $\mathcal{L}^{fr}$ the set of ambient isotopy classes of unoriented framed links (including the empty link $\varnothing$) in $M$, $R$ a commutative ring with unity, and $A \in R$ a fixed invertible element. In addition, let $R\mathcal{L}^{fr}$ be the free $R$-module generated by $\mathcal{L}^{fr}$ and $S_{2, \infty}^{sub}$ the submodule of $R\mathcal{L}^{fr}$ generated by all (local) skein expressions of the form: 
\begin{itemize}
    \item [(i)] $L_+ - AL_0 - A^{-1}L_{\infty}$, and
    \item [(ii)] $L \sqcup \bigcirc  + (A^2 + A^{-2})L$,
\end{itemize} where $\bigcirc$ denotes the trivial framed knot and the {\it skein triple} $(L_+$, $L_0$, $L_{\infty})$ denotes three framed links in $M$, which are identical except in a small $3$-ball in $M$ where they differ as shown in Figure \ref{skeintriple}.

\begin{figure}[h]
    \centering

\[  \begin{minipage}{1.5in} \includegraphics[width=\textwidth]{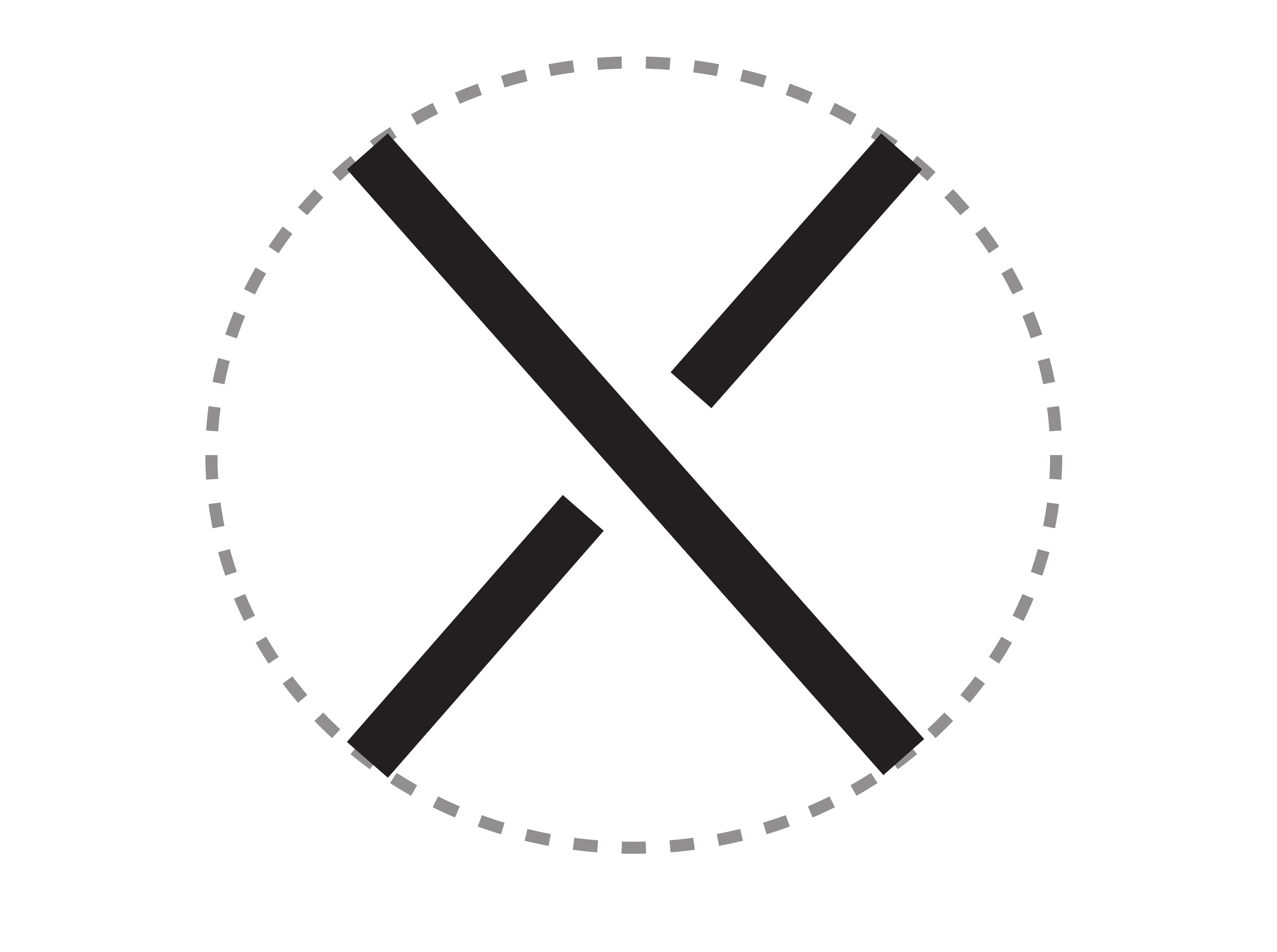} \vspace{-15pt} \[\pmb L_+\] \end{minipage} 
               \qquad
        \begin{minipage}{1.5in}\includegraphics[width=\textwidth]{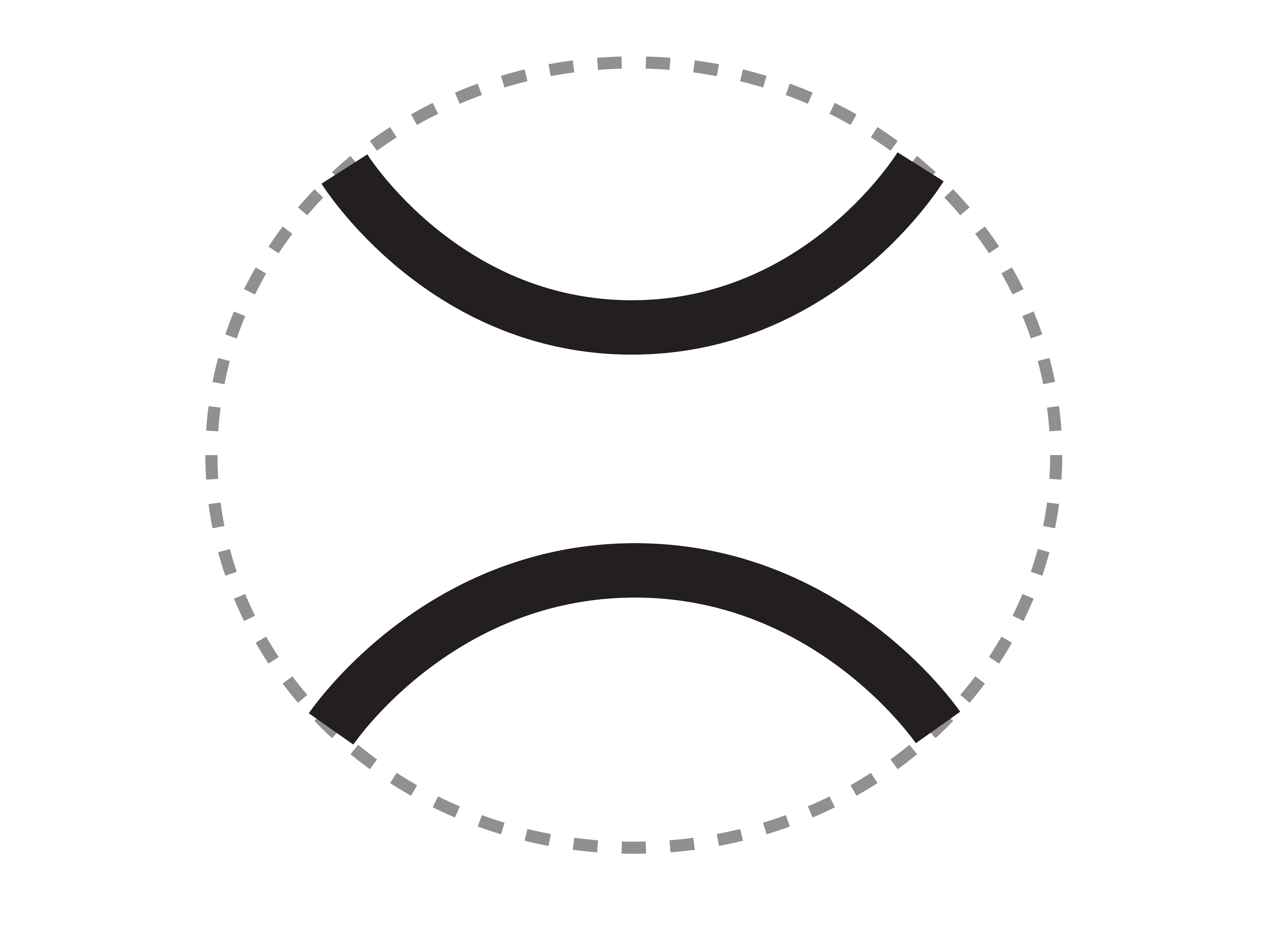} \vspace{-15pt} \[\pmb L_0\] \end{minipage}
         \qquad
        \begin{minipage}{1.5in}\includegraphics[width=\textwidth]{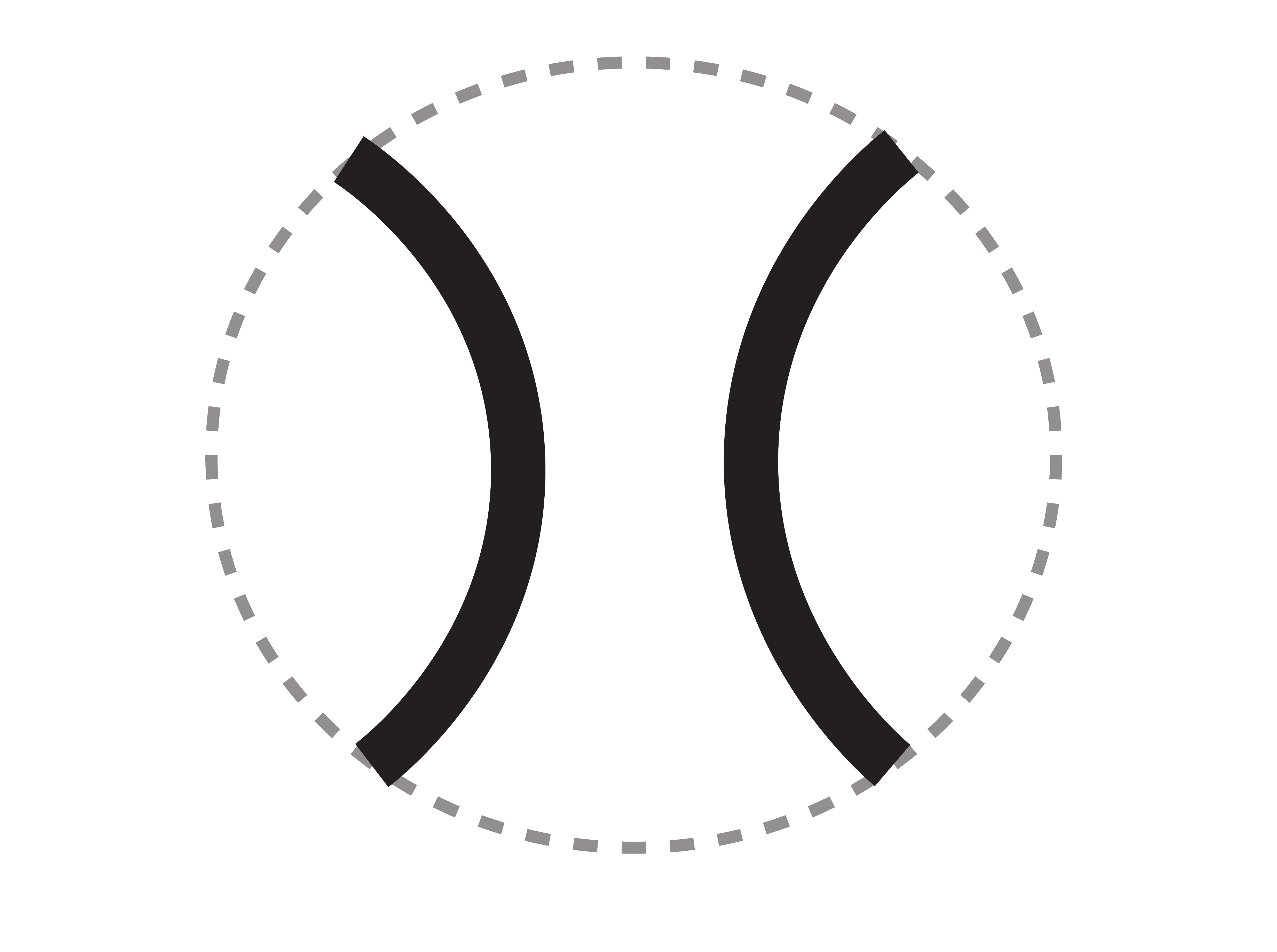} \vspace{-15pt} \[ \pmb L_\infty\]\end{minipage} 
        \]
        \caption{Skein triple for the Kauffman bracket skein module}
          \label{skeintriple}
        \end{figure}
        
The Kauffman bracket skein module (KBSM) of $M$ is defined as the quotient $$\mathcal{S}_{2,\infty}(M;R,A) = R\mathcal{L}^{fr}/S_{2, \infty}^{sub}.$$

\end{definition}

For simplicity, we use the notation $\mathcal{S}_{2,\infty}(M)$ when $R = \mathbb{Z}[A^{\pm 1}]$.

\section{Marché's Conjecture} 


For a long time the structure of the KBSM of $3$-manifolds was elusive and there are few examples for which the structure is completely known. Computations are usually hard and have been carried out over various commutative rings. Over the polynomial ring $\mathbb{Z}[A^{\pm 1}]$, the KBSM has been computed for $S^3$, $I$-bundles over surfaces \cite{smof3,fundamentals}, lens spaces $L(p,q)$ when $p\geq 1$ \cite{kbsmlens}, $S^1 \times S^2$ \cite{s1s2}, the Whitehead manifold \cite{whitehead}, the product of the pair of pants with $S^1$ \cite{pairofpantss1}, certain classes of prism manifolds \cite{prism}, twist knot exteriors \cite{knotext}, the exteriors of torus knots of type $(2,2p+1)$ \cite{comp22p1}, and $3$-manifolds obtained by integral surgery on the right-handed trefoil knot \cite{trefoilbullock}. Over $\mathbb{C}[A^{\pm1}]$, the KBSM of the exteriors of $2$-bridge knots and links was computed in \cite{2bk, 2bl}. When the ring $\mathbb{Z}[A^{\pm 1}]$ is localised by inverting all cyclotomic polynomials, the formula for the KBSM for the quaternionic manifold is known \cite{quaternion}, while over $\mathbb{Q}(A)$, the field of rational functions in the variable $A$, the KBSMs of $T^3$ \cite{min9,li9}, trivial $S^1$-bundles over oriented surfaces with genus $g \geq 2$ \cite{basissurfacetimess1}, and some infinite families of hyperbolic $3$-manifolds \cite{dethyperbolic} are known. Sometime in 2014 or 2015, E. Witten conjectured that the Kauffman bracket skein module over $\mathbb{Q}(A)$ for any closed $3$- manifold is always finite dimensional. This conjecture has never been mentioned in any published work by Witten himself. The first written documentation of it appears in \cite{min9}.  Section $8$ in \cite{surfacetimess1} also briefly discusses this conjecture. Witten's conjecture comes as a surprise in the theory of skein modules since $\mathcal S_{2,\infty}(M, \mathbb C, -1)$ is related to the $SL(2, \mathbb C)$ character variety of $M$, which is usually infinite dimensional.\\
 
In \cite{connsum}, Przytycki showed that if $\mathcal{S}_{2,\infty}(M; \mathbb{Q}(A))$ and $\mathcal{S}_{2,\infty}(N; \mathbb{Q}(A))$ have free parts that are finite dimensional, then so does $\mathcal{S}_{2,\infty}(M \ \# \ N; \mathbb{Q}(A))$, thus showing that Witten's conjecture is stable under connected sums. Using factorisation algebras, the representation theory of quantum groups, and deformation quantization modules, Gunningham, Jordan, and Safronov resolved this conjecture in the affirmative. 

\begin{theorem}\cite{wittenresolved}\label{wittenresolved}

The Kauffman bracket skein module of any closed oriented $3$-manifold over the field $\mathbb{C}(A)$ is finite dimensional. 

\end{theorem}

We note that the KBSM is not the only skein module that exhibits finite dimensionality over the field of rational functions. For closed oriented $3$-manifolds, the $q$-homology skein module also possesses this property. This was proved by Przytycki in \cite{qanalogue}.  
It must also be noted that when $R = \mathbb{Z}[A^{\pm 1}]$, the KBSM of closed oriented $3$-manifolds need not be finite dimensional. The KBSM of $S^1 \times S^2$ is a quintessential example of this phenomenon. In particular, Hoste and Przytycki \cite{s1s2} showed that $\mathcal S_{2,\infty}(S^1 \times S^2)$ has infinitely generated torsion. Thus, much more is known about the structure of the KBSM over $\mathbb Q(A)$ than over $\mathbb Z[A^{\pm 1}]$. March\'e proposed (see \cite{basissurfacetimess1}) the following conjecture about the structure of the KBSM over $\mathbb{Z}[A^{\pm 1}]$. 

\begin{conjecture}\label{marcheconj}

Let $M$ be a closed compact oriented $3$-manifold. Then there exist finitely generated $\mathbb{Z}[A^{\pm 1}]$-modules $N_k$ and an integer $d \geq 0$ such that 
$$\mathcal S_{2, \infty}(M)= (\mathbb{Z}[A^{\pm 1}])^d \oplus \bigoplus \limits_{k \geq 1} N_k,$$ where $N_k$ is an $(A^k - A^{-k})$-torsion module for each $k \geq 1$. 

\end{conjecture}  

This conjecture is true for all the manifolds listed at the beginning of this section. For example, for the $3$-manifold $S^1 \times S^2$, $d = 1$ and $N_k = \cfrac{\mathbb{Z}[A^{\pm 1}]}{1 - A^{2k}}$ for $k\geq 2$. For the $3$-manifold $F \times S^1$, where $F$ is an oriented surface with genus $g\geq 2$, $d = 2^{2g+1}+2g -1$, and as a byproduct of the proof in \cite{basissurfacetimess1}, the torsion elements are always of $(A^k - A^{-k})$-type for some $k \geq 1$.

\section{The Counterexample}

It turns out, however, that Conjecture \ref{marcheconj} is not true in general. While writing \cite{mybook1}, the author discovered a counterexample to this conjecture given by the KBSM of the connected sum of real projective spaces.

\begin{theorem}

Conjecture \ref{marcheconj} is not true when $ M= \mathbb{R}P^3 \ \# \ \mathbb{R}P^3$. In particular, $\mathcal{S}_{2,\infty}(\mathbb{R}P^3 \ \# \ \mathbb{R}P^3)$ does not split into the direct sum of free modules and torsion modules. 

\end{theorem}

We now describe the structure of the KBSM of $\mathbb{R}P^3 \ \# \ \mathbb{R}P^3$. A key ingredient in calculating  $\mathcal{S}_{2,\infty}(\mathbb{R}P^3 \ \# \ \mathbb{R}P^3)$ is the concept of depicting links in $\mathbb{R}P^3 \ \# \ \mathbb{R}P^3$ using arrow diagrams, which we briefly discuss.

\subsection{Arrow diagrams in \pmb {$\mathbb RP^2 \ \hat\times \ S^1$}} In \cite{pairofpantss1}, nice diagrammatic representations of unoriented links in trivial circle bundles over orientable surfaces were introduced by D\c{a}bkowski and Mroczkowski\footnote{A precursor to arrow diagrams was introduced by Turaev in \cite{turaevshadow}.}. These representations are known as arrow diagrams. These arrow diagrams were used to compute the KBSM of the trivial $S^1$-bundle over the pair of pants. In \cite{rp3rp3}, Mroczkowski extended these representations to links in twisted circle bundles over unorientable surfaces (see also \cite{gabmro}). Note that $\mathbb RP^3 \smallsetminus B^3 \cong \mathbb RP^2 \hat \times I$, where $\mathbb RP^2 \hat \times I$ is the twisted $I$-bundle over the real projective plane. Thus, $\mathbb RP^3 \ \# \ \mathbb RP^3 \cong \mathbb RP^2 \ \hat\times \ S^1$, where $\mathbb RP^2 \hat \times S^1$ is the twisted circle bundle over $\mathbb RP^2$. \\ 

\begin{figure}[h]
    \centering
    \begin{minipage}{3in}
    \vspace{0.1pt} \[F \times \{1\}\]\includegraphics[width=\textwidth]{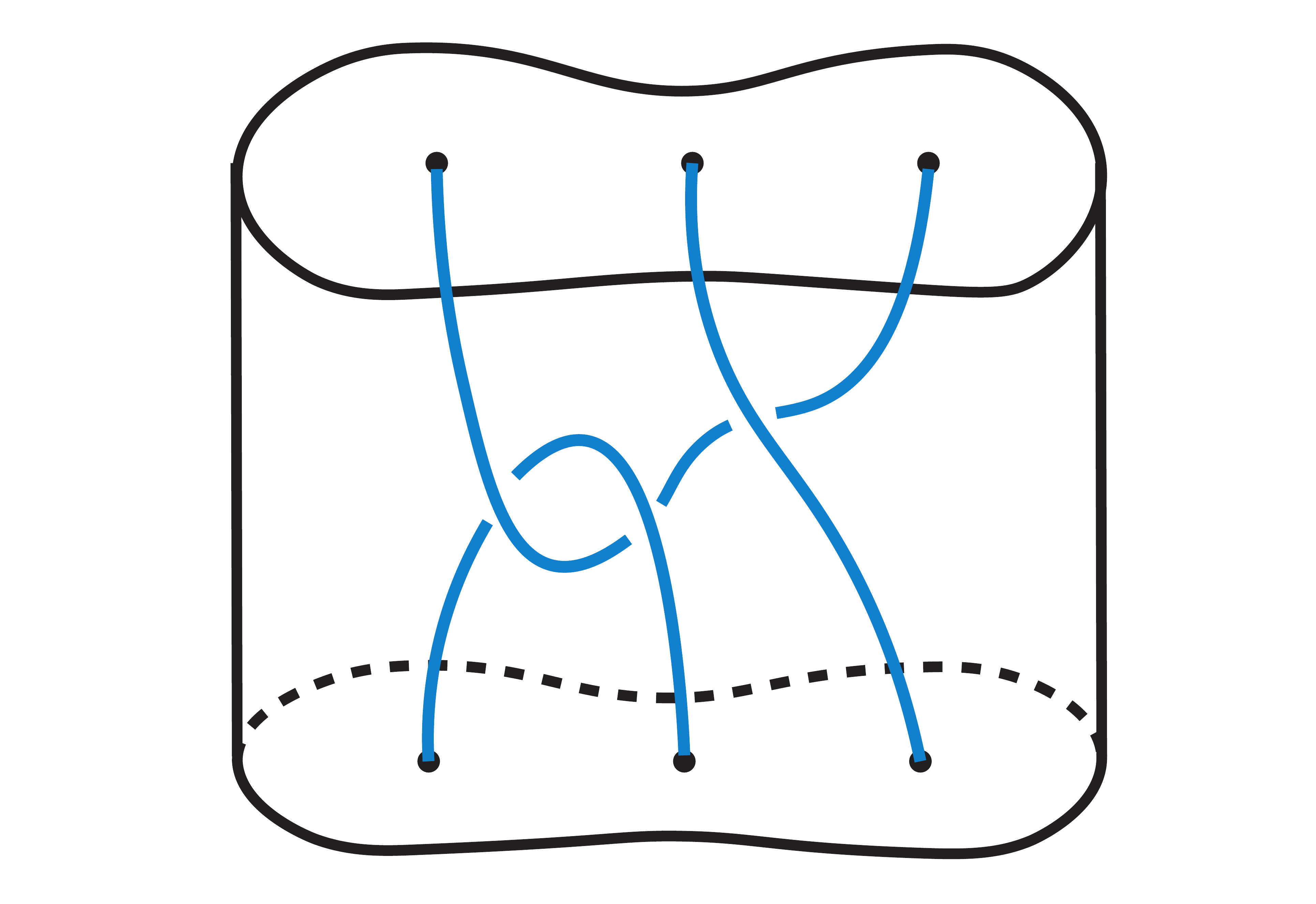} \vspace{-15pt} \[F \times \{0\}\] \end{minipage} $\underrightarrow{\ \ \ \ \ \pi \ \ \ \ \ }$ \ \ \ \ \ \ \   \begin{minipage}{2.5in}\includegraphics[width=\textwidth]{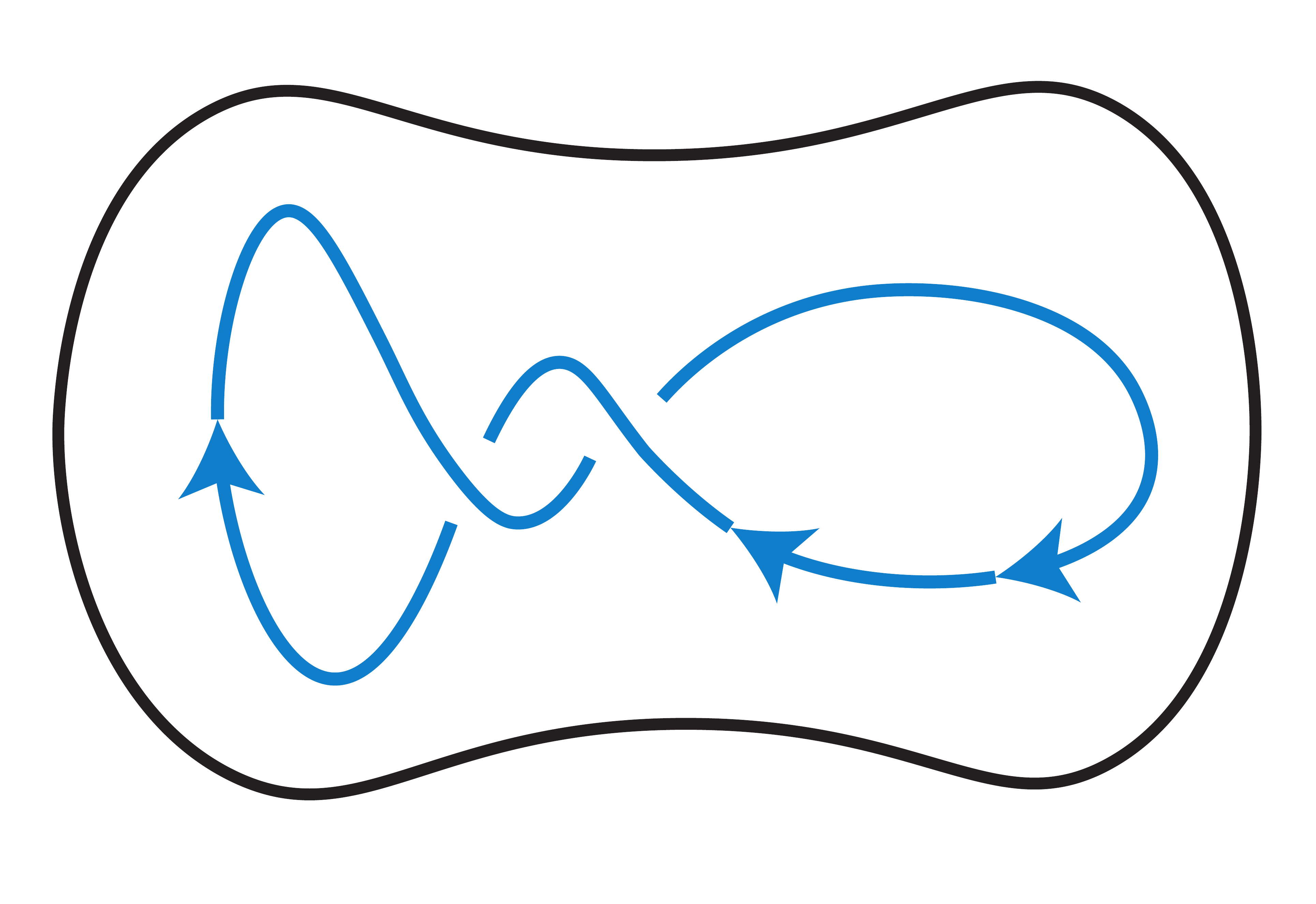} \vspace{-30pt} \[F\] \end{minipage}
    \caption{A link and its arrow diagram in $F \times S^1$}
    \label{arrowdiagramfs1}
\end{figure}

Let $F$ be an orientable surface and consider an unoriented link $L \hookrightarrow F \times S^1$. Cut $F \times S^1$ along $F \times \{1\}$, where $1 \in S^1$. In doing so, we obtain the manifold $F \times [0,1]$ and a collection of circles and arcs that come from $L$. Denote this collection by $L'$. By a general position argument we may assume that $L$ intersects $F \times \{1\}$ transversely (but not orthogonally) in a finite number of points. Consider the projection $\pi : F \times [0,1] \longrightarrow F \times \{0\}$. Projecting the arcs onto $F \times \{0\}$ gives us a set of closed curves in $F$. By another general position argument, we may assume that $\pi(L')$ contains transverse double points and, under $\pi$, the endpoints of the arcs in $L'$ are projected onto points that are distinct from these double points. The closed curves in $\pi(L')$ are equipped with some additional information. 
Besides the usual undercrossing and overcrossing information based on the appearance of the curves in $F \times [0,1]$, we keep track of where the arcs intersect $F \times \{0\}$. This is done by placing at the points of intersection a dot and an arrow pointing in the direction of increasing height in $[0, 1]$ just before and after the cut was made. Thus, at the points that are projections of the endpoints $(x,0)$ and $(x,1)$ of arcs in $L'$, an arrow is placed in the direction in which the height drops by $1$ in $L'$ when the first coordinate crosses $x$. This means that if we travel along $L$ in the direction of the arrow and cross it, we have essentially travelled through the roof, $F \times \{1\}$, and resurfaced at the floor, $F \times \{0\}$. See Figure \ref{arrowdiagramfs1}. \\

\begin{figure}[h]
    \centering
\begin{subfigure}{.49\textwidth}
\centering
\begin{overpic}[unit=1mm, scale = 0.116]{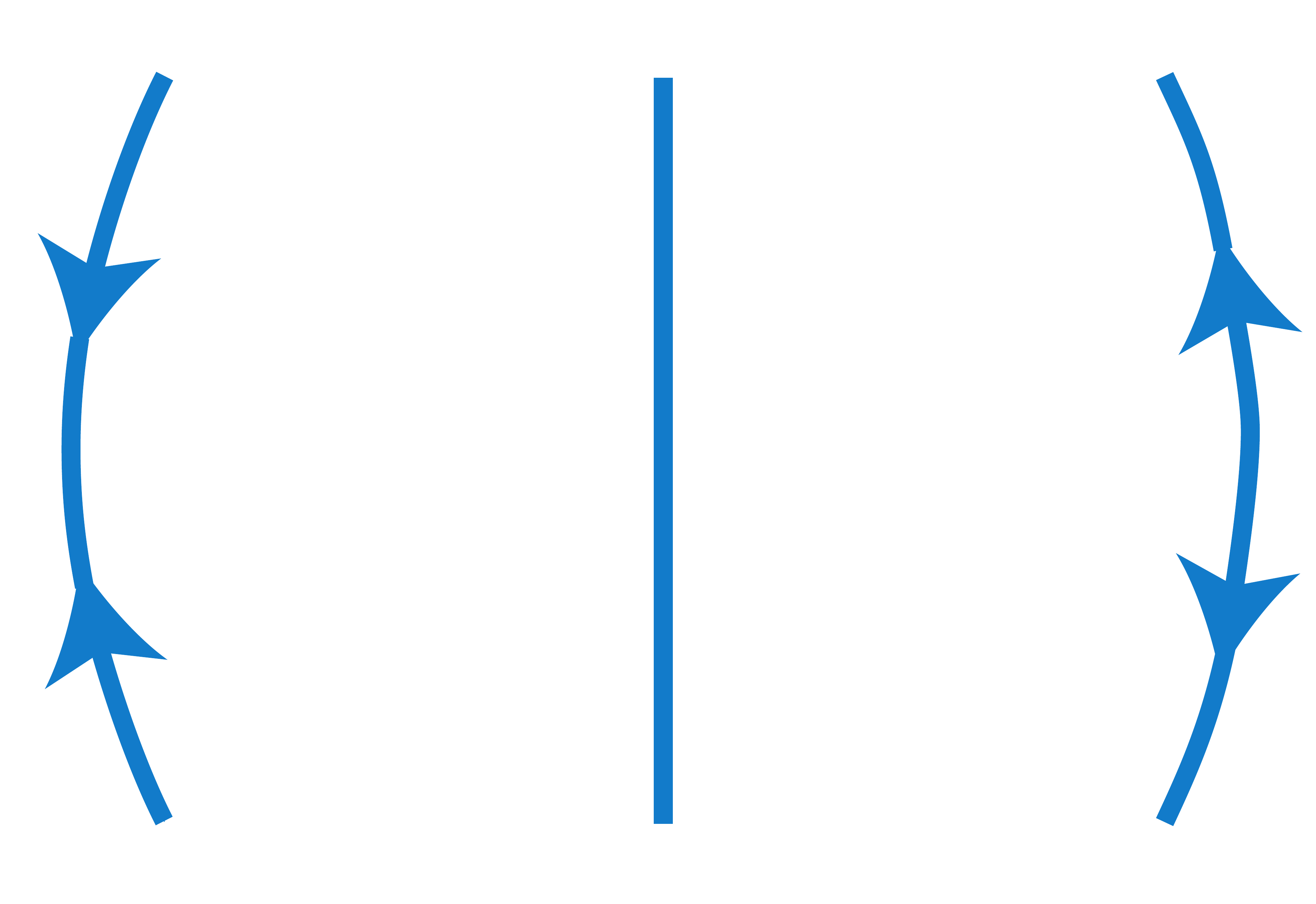}
\put(9,15){$\longleftrightarrow$}
\put(29,15){$\longleftrightarrow$}
\end{overpic} 
\subcaption{}
\end{subfigure}
\centering
\begin{subfigure}{.49\textwidth}
    \begin{minipage}{1.4in}
  \includegraphics[width=\textwidth]{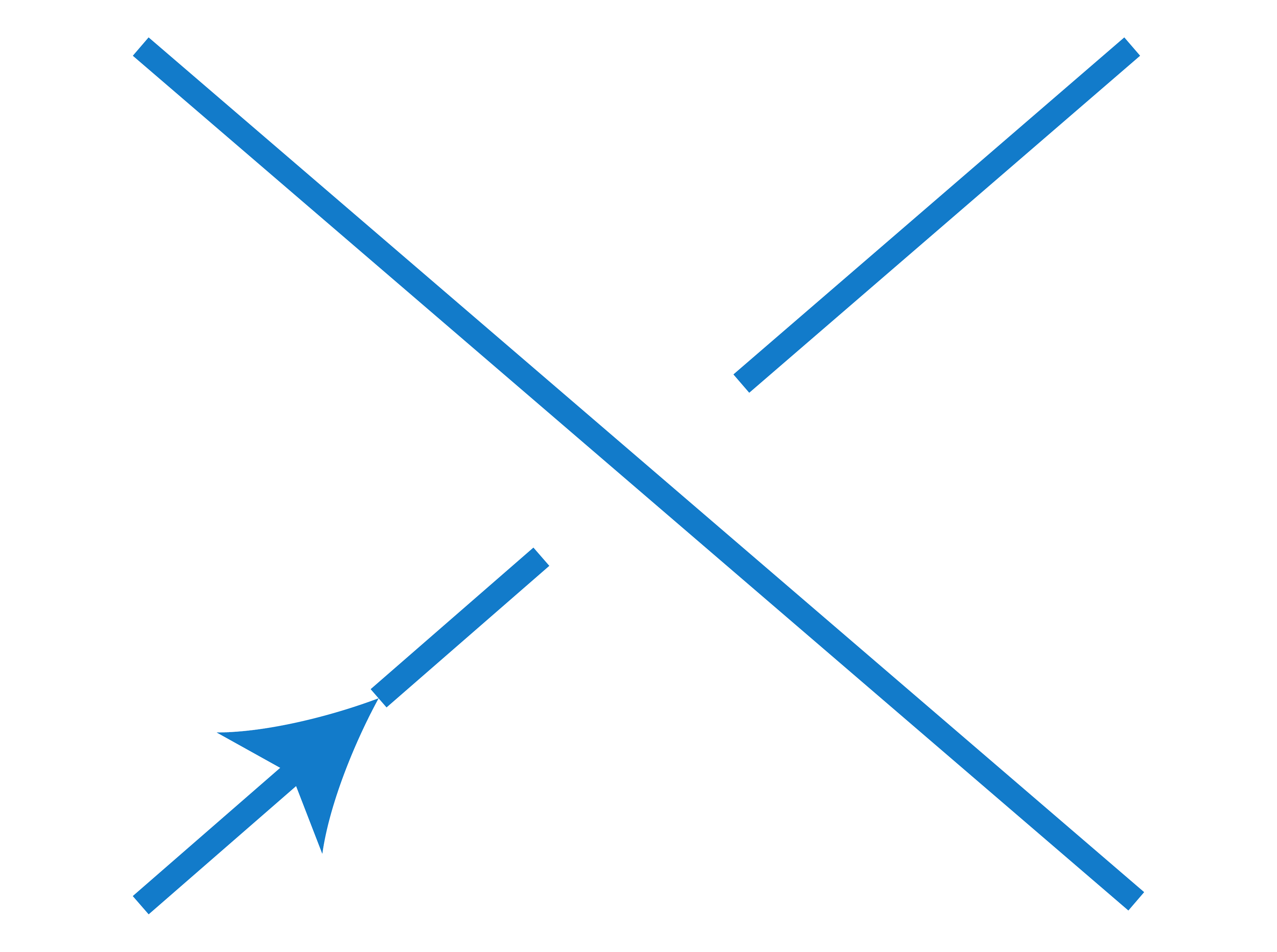}
  \end{minipage}
  $\longleftrightarrow $
  \begin{minipage}{1.4in}
  \includegraphics[width=\textwidth]{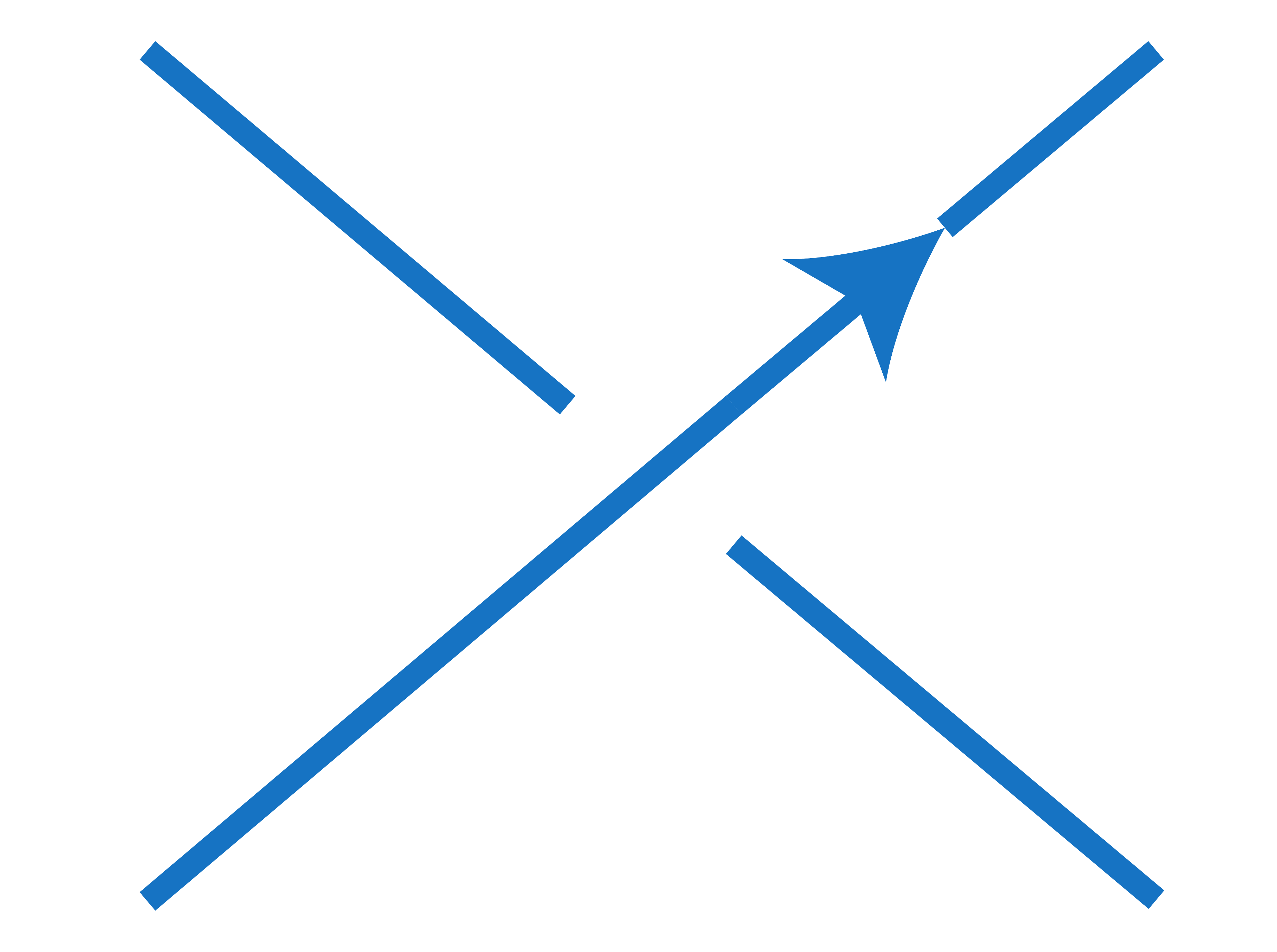} 
  \end{minipage}
  \vspace*{4mm}
  \subcaption{}
\end{subfigure}
\caption{Isotopy moves for arrow diagrams} 
\label{isotopymoves}
    \label{fig:my_label}
\end{figure}

In summary, an arrow diagram $D$ of a link $L$ in $F \times S^1$ is a link diagram in $F$ with crossing information at double points and some arrows. Links in $F \times S^1$ are ambient isotopic to each other if their arrow diagrams are related by a finite sequence of the usual Reidemeister moves and the two additional moves illustrated in Figure \ref{isotopymoves}.  \\

\begin{figure}[h]
    \centering
 \begin{subfigure}{.5\textwidth}
    \begin{minipage}{1.4in}
  \includegraphics[width=\textwidth]{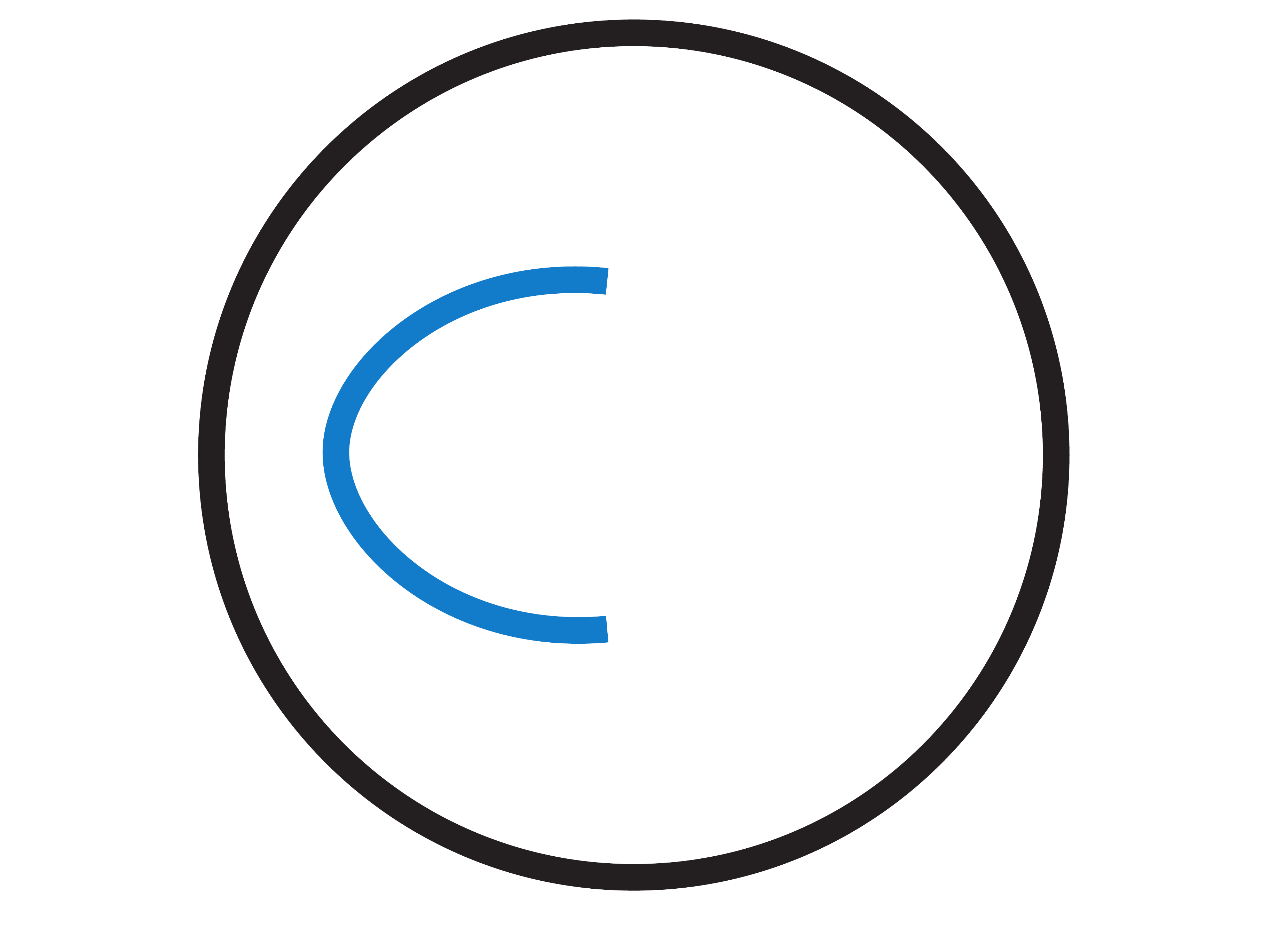}
  \end{minipage}
  $\longleftrightarrow $
  \begin{minipage}{1.4in}
  \includegraphics[width=\textwidth]{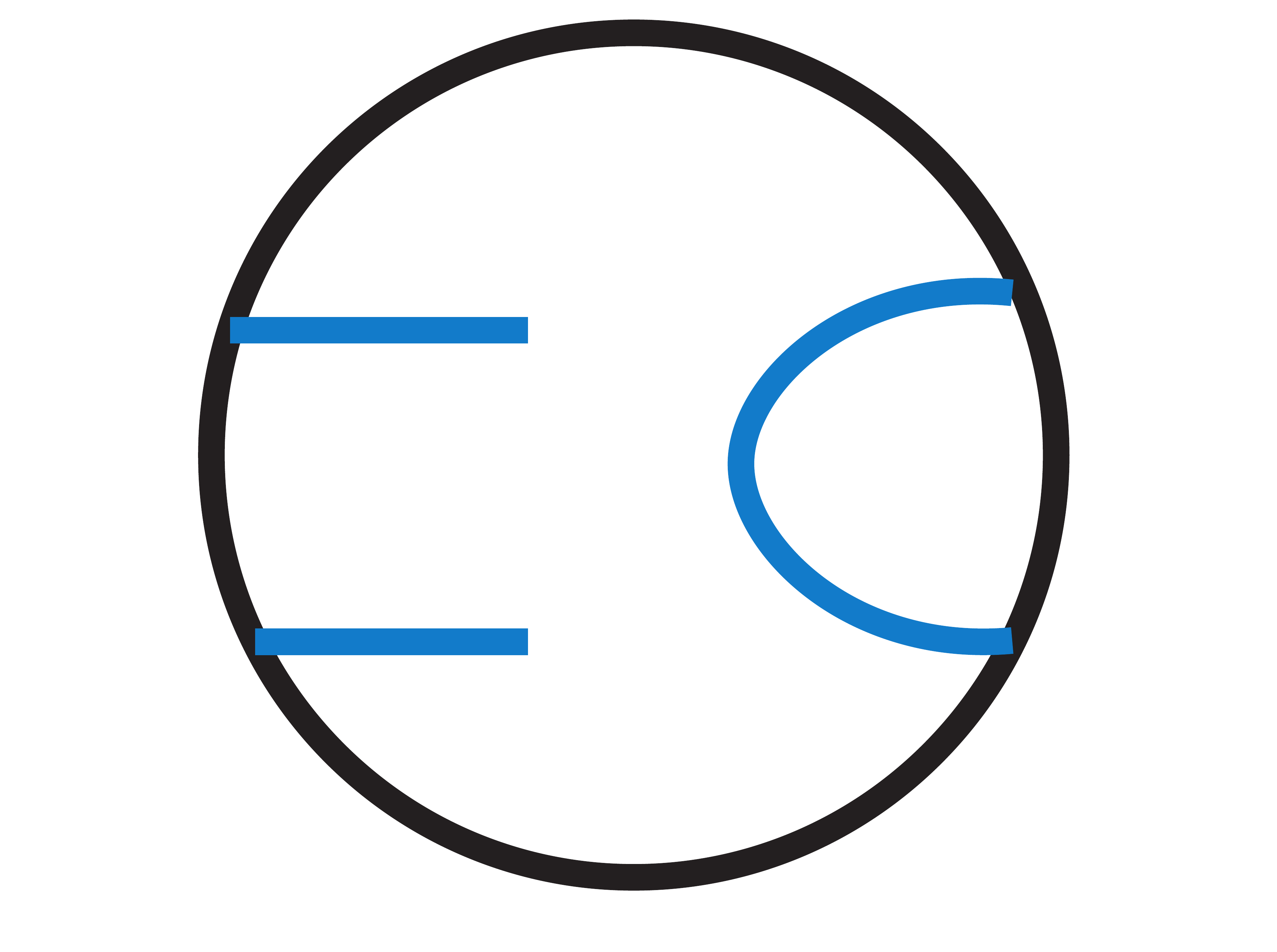} 
  \end{minipage}
  \vspace*{4mm}
  \subcaption{}
\end{subfigure}\hfill
   \centering
\begin{subfigure}{.5\textwidth}
    \begin{minipage}{1.4in}
  \includegraphics[width=\textwidth]{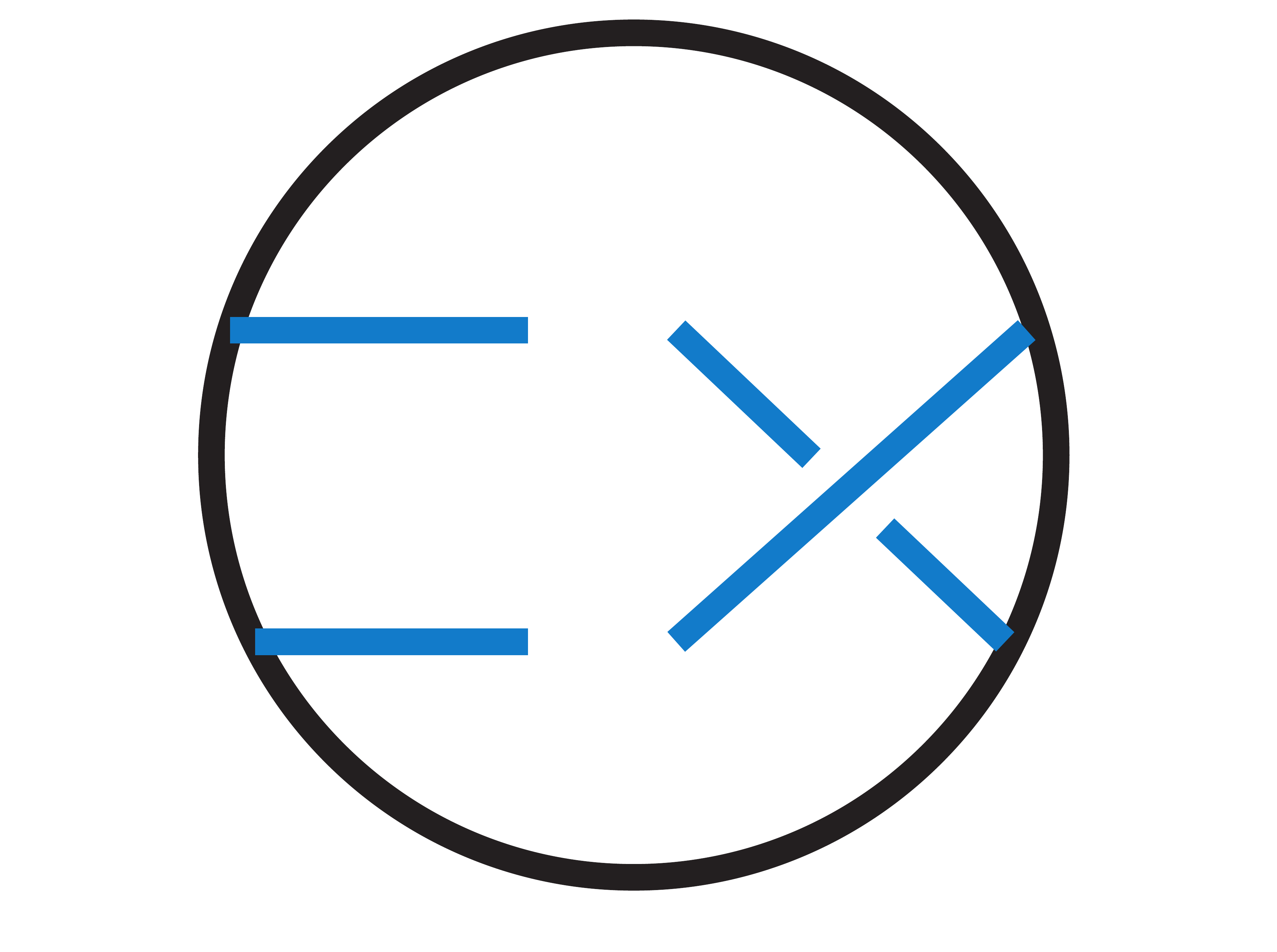}
  \end{minipage}
  $\longleftrightarrow $
  \begin{minipage}{1.4in}
  \includegraphics[width=\textwidth]{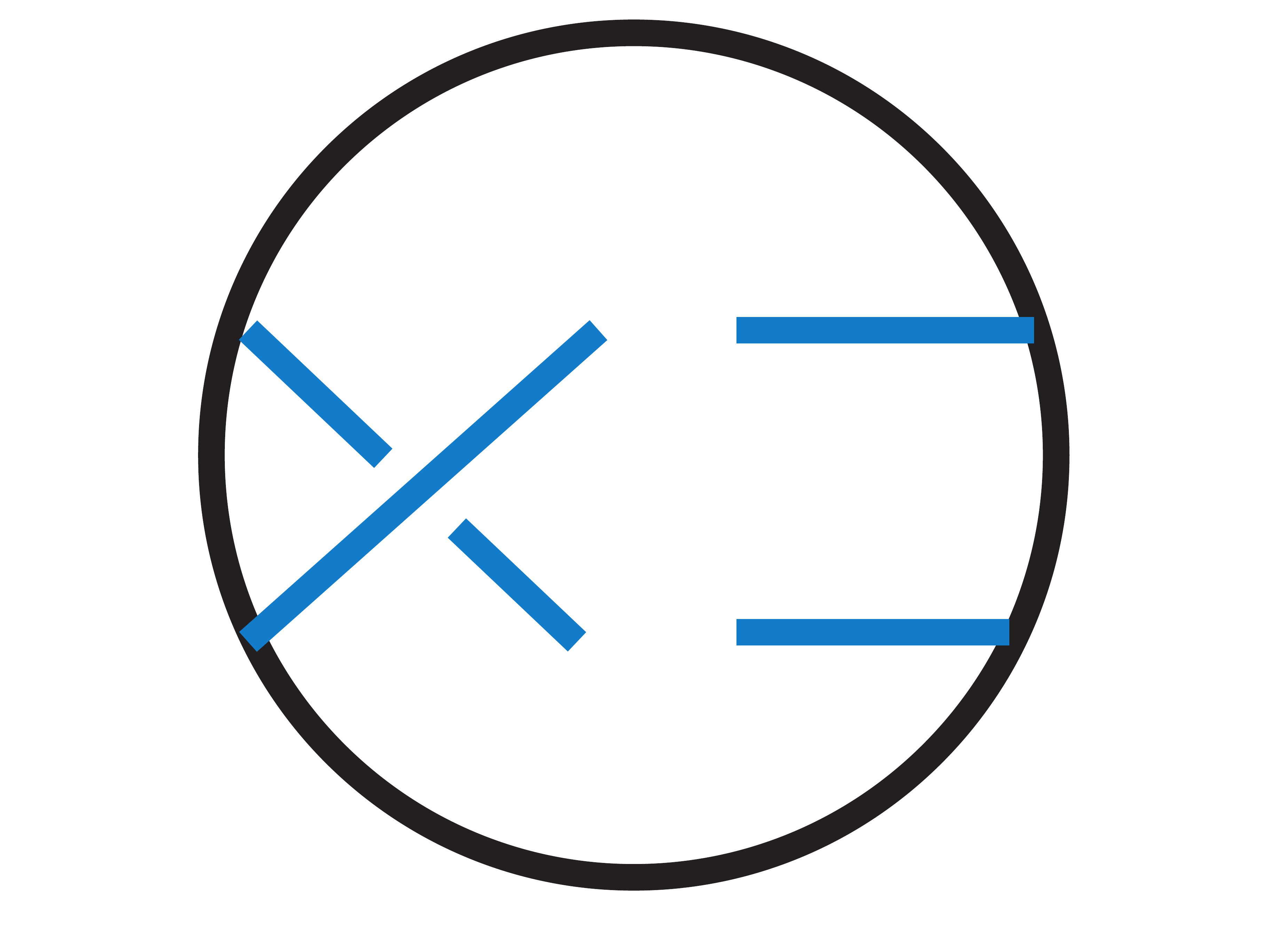} 
  \end{minipage}
  \vspace*{4mm}
  \subcaption{}
\end{subfigure}\hfill
   \centering
\begin{subfigure}{.5\textwidth}
    \begin{minipage}{1.4in}
  \includegraphics[width=\textwidth]{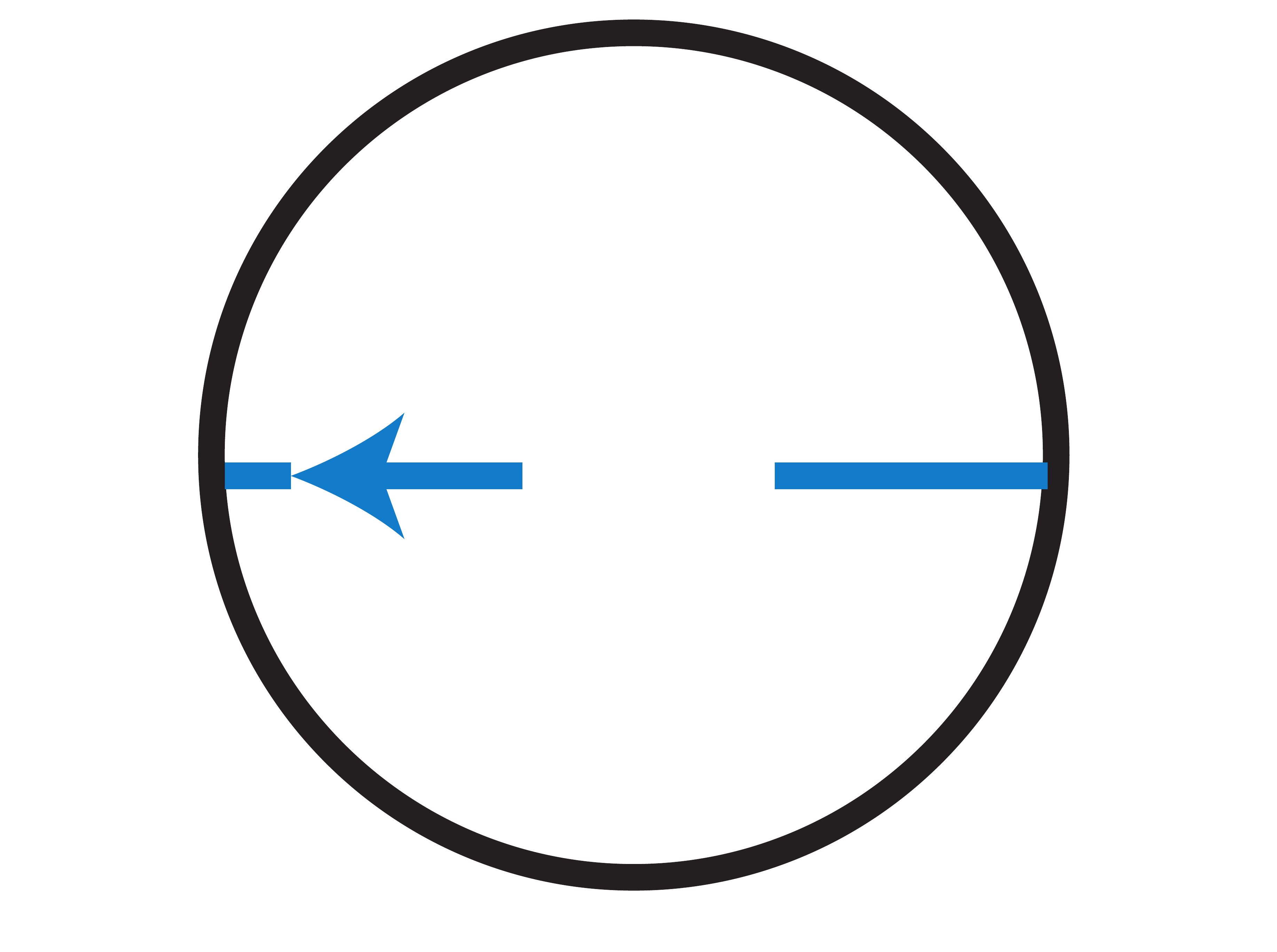}
  \end{minipage}
  $\longleftrightarrow $
  \begin{minipage}{1.4in}
  \includegraphics[width=\textwidth]{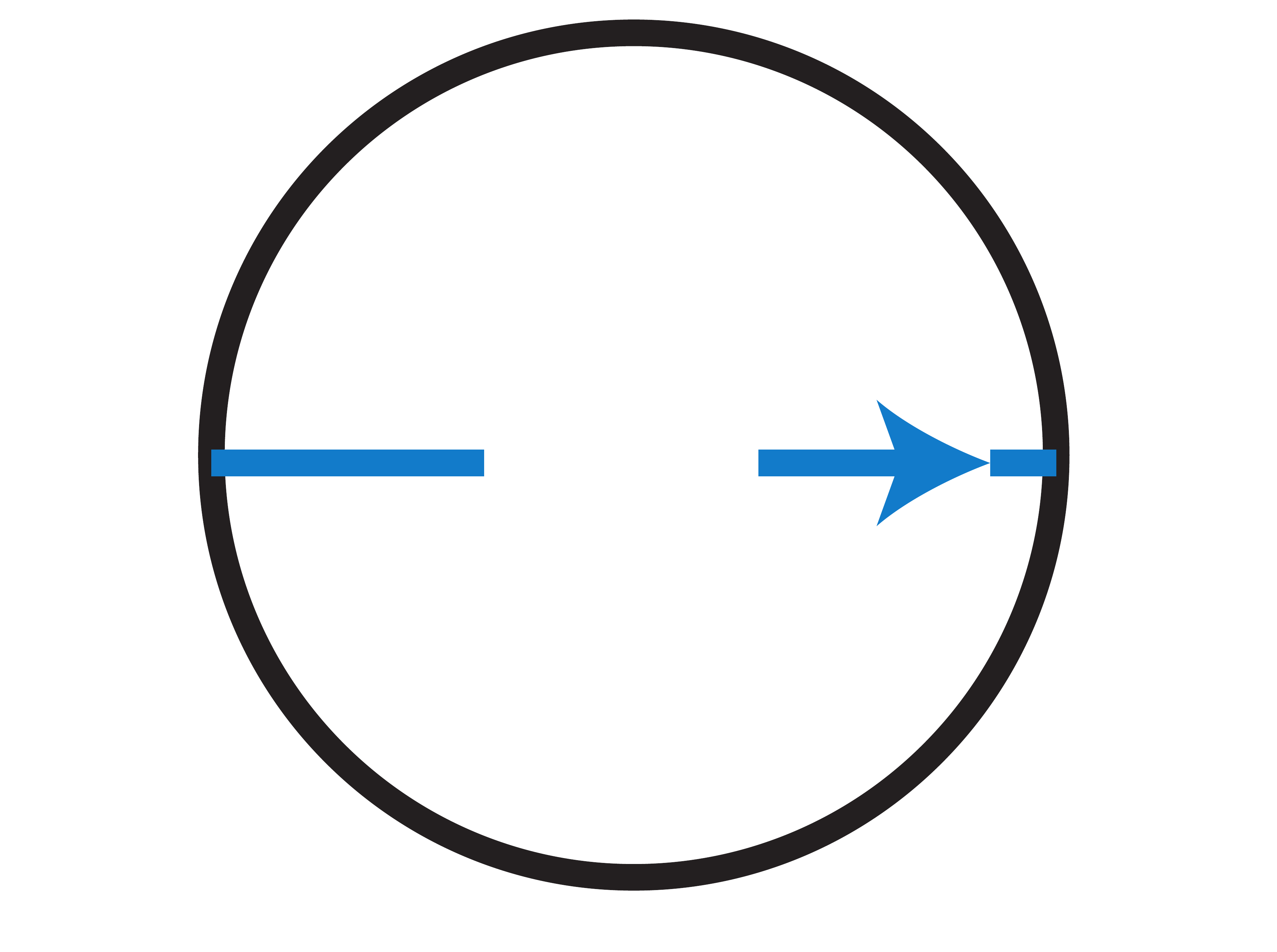} 
  \end{minipage}
  \vspace*{4mm}
  \subcaption{}
\end{subfigure}
\caption{Additional isotopy moves for arrow diagrams in $\mathbb RP^2 \hat \times S^1$}
    \label{hatarrowmoves}
\end{figure}

We can generalise this construction to unoriented links in $\mathbb RP^2 \hat \times S^1$ in the following way. The unorientable surface $\mathbb RP^2$ is obtained by identifying the boundary of $D^2$ using the antipodal map. Thus, $\mathbb RP^2 \hat \times S^1$ is obtained from $D^2 \times S^1$ by identifying $(x,y) \in \partial D^2 \times S^1$ to $(-x,\overline y)$, where $\overline y$ denotes the complex conjugate of $y$. Consider a link $L \hookrightarrow \mathbb RP^2 \hat \times S^1$. In $D^2 \times S^1$, $L$ becomes a collection of closed curves and arcs whose endpoints come in antipodal pairs, $(x,y)$ and $(x, \overline y)$, in $\partial D^2 \times S^1$. As before, we construct arrow diagrams for $L'$ in $D^2 \times S^1$ with the difference that now the arcs have endpoints that appear in antipodal pairs. Thus, arrow diagrams in $\mathbb RP^2 \hat \times S^1$ consist of closed curves and arcs in $D^2$, with the endpoints of the arcs coming in antipodal pairs on $\partial D^2$. Any two links in $\mathbb RP^2 \ \hat\times \ S^1$ are ambient isotopic to each other if their arrow diagrams are related by a finite sequence of the five isotopy moves listed above and the three additional moves illustrated in Figure \ref{hatarrowmoves}. \\

Let $x$ denote the arrow diagram $\vcenter{\hbox{\includegraphics[scale = .4]{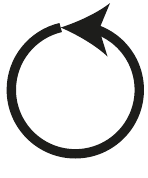}}}$ and $t = -A^{-3}x$, that is, $t$ represents $x$ with a negative full-twist. Then we have the following result about the structure of the KBSM of $\mathbb RP^3 \ \# \ \mathbb RP^3$ over the ring $\mathbb Z[A^{\pm1}]$.

\begin{figure}[h]
    \centering
    \begin{subfigure}{0.35 \textwidth}
    \includegraphics[scale=0.15]{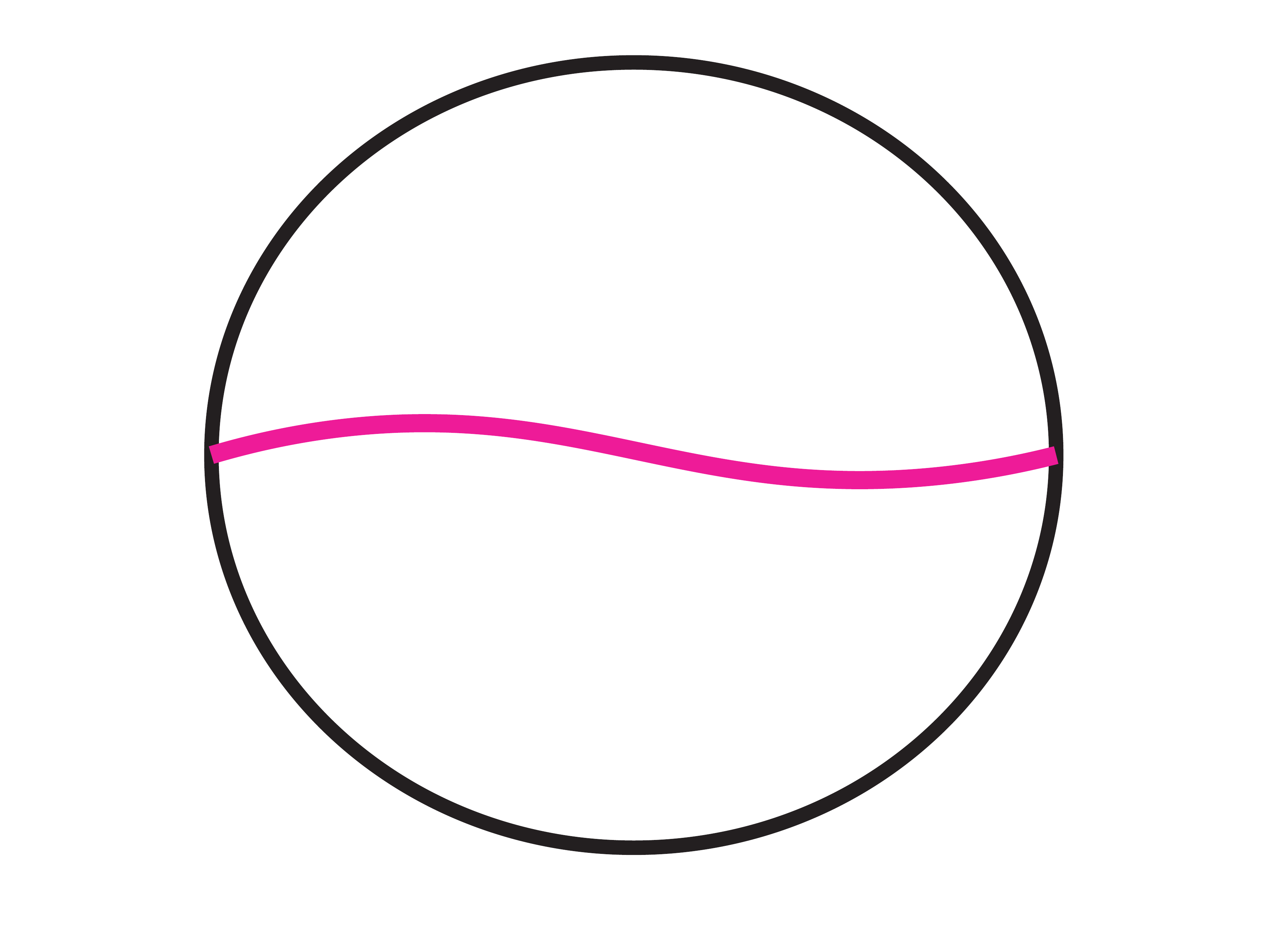}
    \subcaption{$K$}
    \end{subfigure}
     \centering
    \begin{subfigure}{0.35 \textwidth}
    \includegraphics[scale=0.15]{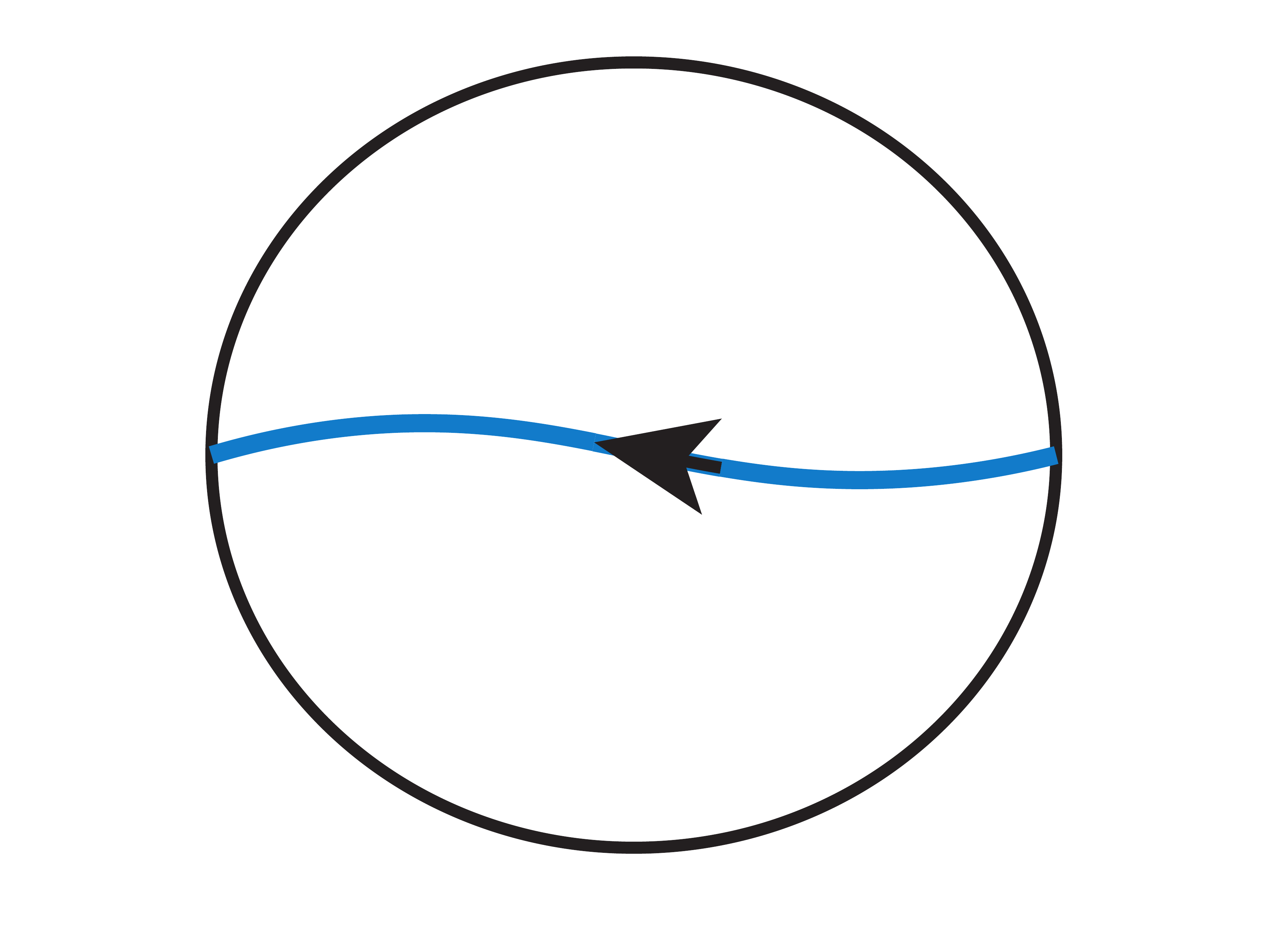}
    \subcaption{$K'$}
    \end{subfigure}
    \caption{The generators of the free part of $\mathcal S_{2,\infty}(\mathbb RP^3 \ \# \ \mathbb RP^3)$}
    \label{freepartgen}
\end{figure}

\begin{theorem}\cite{rp3rp3}\label{rp3rp3}

$\mathcal{S}_{2,\infty}(\mathbb{R}P^3 \ \# \ \mathbb{R}P^3) =\mathbb{Z}[A^{\pm 1}] \oplus \mathbb{Z}[A^{\pm 1}] \oplus \mathbb{Z}[A^{\pm 1}][t]/S$, where $S$ is a submodule of $\mathbb{Z}[A^{\pm 1}][t]$ generated by the following two relations: 
    
    \begin{itemize}
    
    \item [(i)]$(A^{n+1} + A^{-(n + 1)})(S_n(t) - 1)-2(A + A^{-1})\sum \limits_{k=1}^{n/2} A^{n + 2 - 4k}$, for $n \geq 2$ even, 
    
    \item [(ii)] $(A^{n+1} + A^{-(n + 1)})(S_n(t) - t)-2t\sum \limits_{k=1}^{(n-1)/2} A^{n + 1 - 4k}$, for $n \geq 3$ odd.

\end{itemize}

 Here $S_n(t)$ denotes the Chebyshev polynomial of the second kind, which is defined recursively by the equation $S_{n + 1}(t) = t\cdot S_n(t) - S_{n - 1}(t)$, with the initial conditions $S_0(t) = 1$ and $S_1(t) = t$. The generators of the free part are the knots $K$ and $K'$ illustrated in Figure \ref{freepartgen}. $K$ is a knot whose diagram consists of an arc with two antipodal endpoints and no crossings, while $K'$ is a knot whose diagram has an arc with an arrow.

\end{theorem}

   \begin{remark}
    
     $\mathcal{S}_{2,\infty}(\mathbb{R}P^3 \ \# \ \mathbb{R}P^3)$ does not split into the sum of free and torsion parts. See Proposition 4.19 in \cite{rp3rp3}. This is a direct contradiction to Conjecture \ref{marcheconj}, which states that the KBSM of a closed compact oriented $3$-manifold over $\mathbb{Z}[A^{\pm 1}]$ splits into the direct sums of free modules and torsion modules.
    
   \end{remark}
   
       The KBSM of $\mathbb RP^3$ over any commutative ring $R$ with unity has two generators, namely the empty link and the noncontractible curve in $\mathbb RP^2$. An implication of Theorem 1.2 in \cite{connsum} about the KBSM over $\mathbb Q(A)$ of the connected sum of two compact oriented $3$-manifolds is that $\mathcal S_{2,\infty}(\mathbb RP^3 \ \# \ \mathbb RP^3; \mathbb Q(A)) = \mathcal S_{2,\infty}(\mathbb RP^3;\mathbb Q(A)) \otimes \mathcal S_{2,\infty}(\mathbb RP^3;\mathbb Q(A))$. Thus, the dimension of the KBSM of $\mathbb RP^3 \ \# \ \mathbb RP^3$ over $\mathbb Q(A)$ is four. Hence, the rank of the KBSM over $\mathbb Z[A^{\pm 1}]$ of $\mathbb RP^3 \ \# \ \mathbb RP^3$ is also four. Moreover, $\mathcal S_{2,\infty}(\mathbb RP^3 \ \# \ \mathbb RP^3)$ contains torsion elements as seen in \cite{rp3rp3}. Compare with Theorem 4.4 in \cite{algtop}.\\

    \section{Future Directions}
    
    While Marché's conjecture is true for all of the examples listed in Section \ref{intromarche}, the $3$-manifold $\mathbb{R}P^3 \ \# \ \mathbb{R}P^3$ is an outlier. Thus, to understand the structure of the Kauffman bracket skein module of an oriented $3$-manifold over the ring $\mathbb{Z}[A^{\pm 1}]$, one should study that of the connected sums of oriented $3$-manifolds. The first step in doing this would be determining the structure of the Kauffman bracket skein module of the connected sums of handlebodies. In \cite{connsum}, a theorem about the structure of the KBSM of the connected sums of two handlebodies $H
    _n \ \# \ H_m$, $n,m \in \mathbb N$, was stated. A counterexample to the theorem was provided in \cite{counterhandle} for $n+m\geq 3$. In ongoing work, joint with Lê and Przytycki, we proved the following result about the KBSM of the connected sum of two solid tori, that is, when $n=m=1$.
    
\begin{figure}[h]
    \centering
\begin{subfigure}{.5\textwidth}
\begin{overpic}[unit=1mm, scale = 0.5]{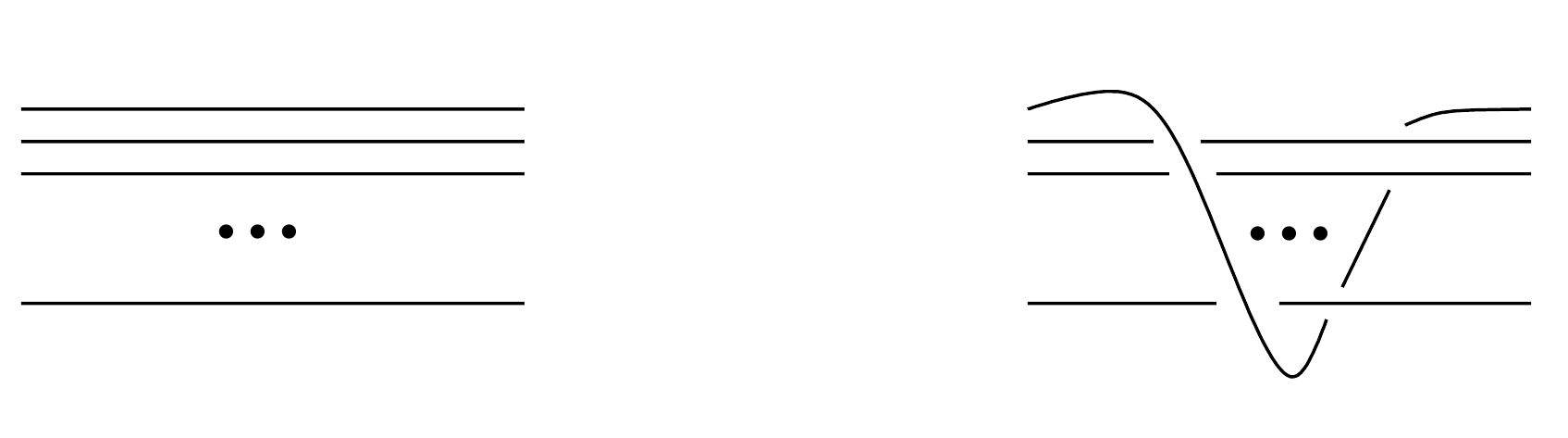}

\put(-8,10){$k$}
\put(-5,11){$\Bigg\{$}
\put(12.5,-1){$z_k$}
\put(67,-1){$u(z_k)$}
\end{overpic}
    \caption{{\color{white}.}}
    \label{fig7.1}
\end{subfigure}
\begin{subfigure}{.5\textwidth}
    \centering
    \begin{overpic}[unit=1mm, scale = 0.52]{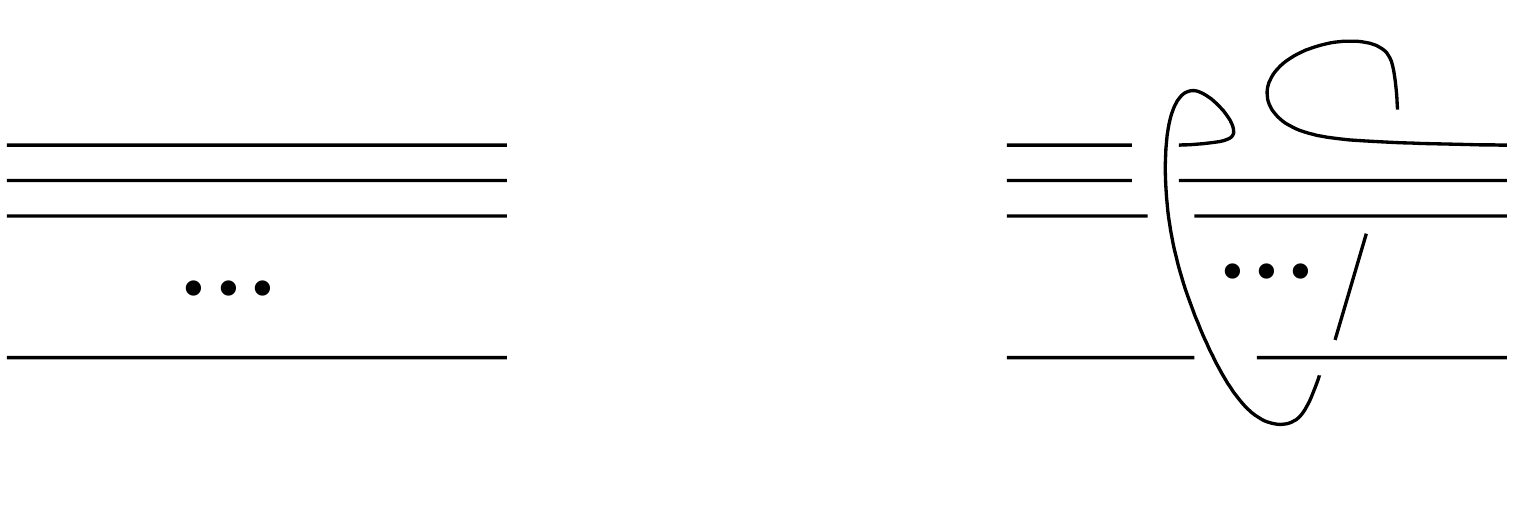}
\put(11,-1){$z_k$}
\put(62,-1){$sl_{\partial D}(z_k)$}
\end{overpic}
    \caption{{\color{white}.}}
    \label{fig7.2}
    \end{subfigure}
\caption{{Illustration of the handle sliding relation that generates $\mathcal I$}}
\end{figure}
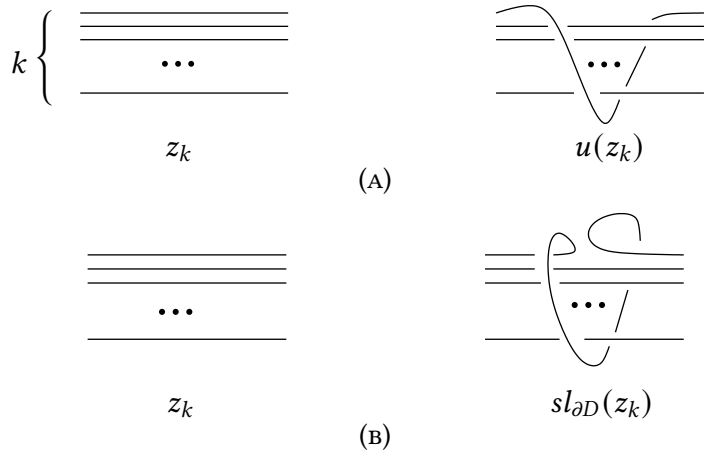

    \begin{theorem}\cite{blp}\label{przytyckitheorem}

Let $H_n$ denote a genus $n$ handlebody and $F_{0,n+1}$ be a disc with $n$ holes so that $H_n = F_{0,n+1} \times I$. Then,
$${\mathcal S_{2,\infty}}(H_1 \ \# \ H_1) = {\mathcal S_{2,\infty}}(H_{2})/\mathcal I.$$
More precisely, the natural epimorphism $i_*: \mathcal{S}_{2,\infty}(H_{2})/\mathcal{I} \longrightarrow \mathcal{S}_{2,\infty}(H_1 \ \# \ H_1)$ is an isomorphism, where $\mathcal I$ is the submodule generated by the expressions $z_k-A^6u(z_k)$, 
for any even $k\geq 2$. Here $z_k \in B_k(F_{0,3})$, where $B_k(F_{0,3})$ is a subset of a basis of $\mathcal{S}_{2,\infty}(F_{0,3} \times I)$ composed of links that have no contractible components and
have geometric intersection number $k$ with a disc $D$ that separates
$H_n$ and $H_m$. $u(z_k)$ is a modification of $z_k$
in the neighbourhood of $D$, as shown in Figure \ref{fig7.1}.
The relation $z_k = A^6u(z_k)$, is a result of the sliding relation
$z_k = sl_{\partial D}(z_k)$ illustrated in Figure \ref{fig7.2}.

\end{theorem}

In future work, the author plans to work on generalising this result to the KBSMs  of $L(p,q) \ \# \ H_1$ and $L(p,q) \ \# \ L(p',q')$, for lens spaces $L(p,q)$ and $L(p',q')$. \\

We also note that the counterexample provided to Marché's conjecture involves the connected sum of oriented $3$-manifolds. We conjecture the following, where part (\ref{conj1}) is stronger than (\ref{conj2}) and (\ref{conj2}) is in turn stronger than (\ref{conj3}): 

\begin{conjecture} \hfill

\begin{enumerate} 
    \item \label{conj1} The KBSM of any closed, prime, oriented $3$-manifold can be decomposed into the direct sum of free modules and torsion modules.
    
    \item \label{conj2} The KBSM of any closed, irreducible, atoroidal, oriented $3$-manifold can be decomposed into the direct sum of free modules and torsion modules. 
    
    \item \label{conj3} The KBSM of any closed, oriented, irreducible, non-Haken\footnote{A non-Haken closed, oriented $3$-manifold is an irreducible manifold that contains no embedded oriented, incompressible surfaces.} $3$-manifold does not contain torsion and is free. This conjecture follows from a conjecture of Przytycki in \cite{kirbyproblemlist} and from Theorem \ref{wittenresolved}.
\end{enumerate}

\end{conjecture}

\section{Acknowledgements}

The author was supported by Dr. Max Rössler, the Walter Haefner Foundation, and the ETH Zürich Foundation.

\end{document}